\documentclass[10.5pt]{amsart}
\usepackage{amsmath}
\usepackage{amssymb}
\usepackage{amscd}
\usepackage{xspace}
\usepackage{verbatim}
\begin{document}
%\numberwithin{equation}{section}
%%%%%%%%%%%%
%%Setting up the TITLE and AUTHOR
\title[Capacitary representation]{Capacitary representation of positive solutions of
semilinear parabolic equations}
\author{Moshe Marcus }
\address{Department of Mathematics, Technion\\
 Haifa, ISRAEL}
%%%%%%%
\author{ Laurent Veron}
\address{Laboratoire de Math\'ematiques, Facult\'e des Sciences\\
Parc de Grandmont, 37200 Tours, FRANCE}
\date{}
%% FONT commands
\newcommand{\txt}[1]{\;\text{ #1 }\;}%% Used in math only
\newcommand{\tbf}{\textbf}%% Bold face. Usage: \tbf{...}
\newcommand{\tit}{\textit}%% Italic
\newcommand{\tsc}{\textsc}%% Small caps
\newcommand{\trm}{\textrm}
\newcommand{\mbf}{\mathbf}%% Math bold
\newcommand{\mrm}{\mathrm}%% Math Roman
\newcommand{\bsym}{\boldsymbol}%% Bold math symbol
%%Macros for changing font size in math.
\newcommand{\scs}{\scriptstyle}%% as in subscript
\newcommand{\sss}{\scriptscriptstyle}%% as in sub-subscript
\newcommand{\txts}{\textstyle}
\newcommand{\dsps}{\displaystyle}
%%Macros for changing font size in text.
\newcommand{\fnz}{\footnotesize}
\newcommand{\scz}{\scriptsize}
%%\tiny<\scz<\fsz<\small<\large<\Large<\huge<\Huge
%%%%%%%%%%%%
%%%%%%%%%%%%
%% EQUATION commands
\newcommand{\be}{\begin{equation}}
\newcommand{\bel}[1]{\begin{equation}\label{#1}}
\newcommand{\ee}{\end{equation}}%% This macro does not work with amstex.
\newcommand{\eqnl}[2]{\begin{equation}\label{#1}{#2}\end{equation}}%%use not
%advisable; confusing
%%%%%%%%%%%%%%%
%% Unnumbered THEOREM env.
%% New env. to be used for unnumbered theorem, lemma etc. (but with specified name)
\newtheorem{subn}{\name}
\renewcommand{\thesubn}{}
\newcommand{\bsn}[1]{\def\name{#1}\begin{subn}}
\newcommand{\esn}{\end{subn}}
%%%%%%%%%%%%%%
%% NUMBERED THEOREM env.
%% Environments: theorem, lemma, corollary defintion and related commands,
%% designed to provide consecutive numbering of these forms.
\newtheorem{sub}{\name}%[section]
\newcommand{\dn}[1]{\def\name{#1}}   %used in conjuction with sub or subn.
\newcommand{\bs}{\begin{sub}}
\newcommand{\es}{\end{sub}}
\newcommand{\bsl}[1]{\begin{sub}\label{#1}}
%% the above must be preceeded by \dn (name definition),
%% however this is superceded by the list of commands bth etc.  below.
%%%%%%%%%%%%
%% NUMBERED THEOREM env. (cont.)
%% List of commands derived from 'sub' env. for theorem, lemma etc.
%% designed to provide consecutive numbering of these forms.
\newcommand{\bth}[1]{\def\name{Theorem}\begin{sub}\label{t:#1}}
\newcommand{\blemma}[1]{\def\name{Lemma}\begin{sub}\label{l:#1}}
\newcommand{\bcor}[1]{\def\name{Corollary}\begin{sub}\label{c:#1}}
\newcommand{\bdef}[1]{\def\name{Definition}\begin{sub}\label{d:#1}}
\newcommand{\bprop}[1]{\def\name{Proposition}\begin{sub}\label{p:#1}}
%%%%%%%%%%%%%%%%%%%%%%%%%%%%%%%%%%
%% RERERENCE commands.
%% \newcommand{\R}[1]{(\ref{#1})}
\newcommand{\R}{\eqref}
\newcommand{\rth}[1]{Theorem~\ref{t:#1}}
\newcommand{\rlemma}[1]{Lemma~\ref{l:#1}}
\newcommand{\rcor}[1]{Corollary~\ref{c:#1}}
\newcommand{\rdef}[1]{Definition~\ref{d:#1}}
\newcommand{\rprop}[1]{Proposition~\ref{p:#1}}
%%%%%%%%%%%
%% ARRAY commands.
\newcommand{\BA}{\begin{array}}
\newcommand{\EA}{\end{array}}
\newcommand{\BAN}{\renewcommand{\arraystretch}{1.2}
\setlength{\arraycolsep}{2pt}\begin{array}}
\newcommand{\BAV}[2]{\renewcommand{\arraystretch}{#1}
\setlength{\arraycolsep}{#2}\begin{array}}
%Note: The first variable gives the amount of stretching: (#1) x default.
%For instance #1=1.2 means a 20% stretching. The second variable should be
%written for instance in the form  4pt ; here the default is 5pt
%\newcommand{\EAN}{\end{array}\setlength{\arraycolsep}{5pt}}
\newcommand{\BSA}{\begin{subarray}}
\newcommand{\ESA}{\end{subarray}}
%Note: These are used in subscripts as well as superscripts. They work essentially
%% like 'array'.
\newcommand{\BAL}{\begin{aligned}}
\newcommand{\EAL}{\end{aligned}}
\newcommand{\BALG}{\begin{alignat}}
\newcommand{\EALG}{\end{alignat}}%% the abbrev. does not work with latex2e
\newcommand{\BALGN}{\begin{alignat*}}
\newcommand{\EALGN}{\end{alignat*}}%% the abbrev. does not work with latex2e
%%%%%%%%%
%% PROOF, REMARF etc.
\newcommand{\note}[1]{\textit{#1.}\hspace{2mm}}
\newcommand{\Remark}{\note{Remark}}
%%%%%%%% Style command.
\newcommand{\modin}{$\,$\\[-4mm] \indent}
%% To be used after \mysection in order to start new line with \indent.
%%%%%%%%%%%%
%% MATHEMATICAL symbols
\newcommand{\forevery}{\quad \forall}
\newcommand{\set}[1]{\{#1\}}
\newcommand{\setdef}[2]{\{\,#1:\,#2\,\}}
\newcommand{\setm}[2]{\{\,#1\mid #2\,\}}
%% Arrows
\newcommand{\lra}{\longrightarrow}
\newcommand{\lla}{\longleftarrow}
\newcommand{\llra}{\longleftrightarrow}
\newcommand{\Lra}{\Longrightarrow}
\newcommand{\Lla}{\Longleftarrow}
\newcommand{\Llra}{\Longleftrightarrow}
\newcommand{\warrow}{\rightharpoonup}
%% Brackets, delimiters
\newcommand{\paran}[1]{\left (#1 \right )}%% adjustable parantheses
\newcommand{\sqbr}[1]{\left [#1 \right ]}%% adjustable square brackets
\newcommand{\curlybr}[1]{\left \{#1 \right \}}%% adjustable curly brackets
\newcommand{\abs}[1]{\left |#1\right |}%% adjustable vertical delimiters
\newcommand{\norm}[1]{\left \|#1\right \|}%% adjustable norm
\newcommand{\paranb}[1]{\big (#1 \big )}%% non-adjustable parantheses (big)
\newcommand{\lsqbrb}[1]{\big [#1 \big ]}%% non-adjustable square brackets (big)
\newcommand{\lcurlybrb}[1]{\big \{#1 \big \}}%% non-adjustable curly brackets (big)
\newcommand{\absb}[1]{\big |#1\big |}%% non-adjustable vertical delimiters (big)
\newcommand{\normb}[1]{\big \|#1\big \|}%% non-adjustable norm (big)
\newcommand{\paranB}[1]{\Big (#1 \Big )}%% non-adjustable parantheses (Big)
\newcommand{\absB}[1]{\Big |#1\Big |}%% non-adjustable vertical delimiters (Big)
\newcommand{\normB}[1]{\Big \|#1\Big \|}%% non-adjustable norm (Big)

%% More mathematical symbols
\newcommand{\thkl}{\rule[-.5mm]{.3mm}{3mm}}
\newcommand{\thknorm}[1]{\thkl #1 \thkl\,}
\newcommand{\trinorm}[1]{|\!|\!| #1 |\!|\!|\,}
\newcommand{\bang}[1]{\langle #1 \rangle}%% angle bracket
\newcommand{\vstrut}[1]{\rule{0mm}{#1}}
\newcommand{\rec}[1]{\frac{1}{#1}}
%% OPERATOR names.
\newcommand{\opname}[1]{\mbox{\rm #1}\,}
\newcommand{\supp}{\opname{supp}}
\newcommand{\dist}{\opname{dist}}
\newcommand{\myfrac}[2]{{\displaystyle \frac{#1}{#2} }}
\newcommand{\myint}[2]{{\displaystyle \int_{#1}^{#2}}}
%%%%%%%%%%
%%%%%%% SPACE commands
\newcommand{\q}{\quad}
\newcommand{\qq}{\qquad}
\newcommand{\hsp}[1]{\hspace{#1mm}}
\newcommand{\vsp}[1]{\vspace{#1mm}}
%%%%%%%%%%%
%% ABREVIATIONS
\newcommand{\ity}{\infty}
\newcommand{\prt}{\partial}
\newcommand{\sms}{\setminus}
\newcommand{\ems}{\emptyset}
\newcommand{\ti}{\times}
\newcommand{\pr}{^\prime}
\newcommand{\ppr}{^{\prime\prime}}
\newcommand{\tl}{\tilde}
\newcommand{\sbs}{\subset}
\newcommand{\sbeq}{\subseteq}
\newcommand{\nind}{\noindent}
\newcommand{\ind}{\indent}
\newcommand{\ovl}{\overline}
\newcommand{\unl}{\underline}
\newcommand{\nin}{\not\in}
\newcommand{\pfrac}[2]{\genfrac{(}{)}{}{}{#1}{#2}}% frac with parantheses.
%%%%%%%%%%%
%%%%%%%%%%%%%

%%Macros for Greek letters.
\def\ga{\alpha}     \def\gb{\beta}       \def\gg{\gamma}
\def\gc{\chi}       \def\gd{\delta}      \def\ge{\epsilon}
\def\gth{\theta}                         \def\vge{\varepsilon}
\def\vgf{\phi}       \def\vgf{\varphi}    \def\gh{\eta}
\def\gi{\iota}      \def\gk{\kappa}      \def\gl{\lambda}
\def\gm{\mu}        \def\gn{\nu}         \def\gp{\pi}
\def\vgp{\varpi}    \def\gr{\rho}        \def\vgr{\varrho}
\def\gs{\sigma}     \def\vgs{\varsigma}  \def\gt{\tau}
\def\gu{\upsilon}   \def\gv{\vartheta}   \def\gw{\omega}
\def\gx{\xi}        \def\gy{\psi}        \def\gz{\zeta}
\def\Gg{\Gamma}     \def\Gd{\Delta}      \def\vgf{\Phi}
\def\Gth{\Theta}
\def\Gl{\Lambda}    \def\Gs{\Sigma}      \def\Gp{\Pi}
\def\Gw{\Omega}     \def\Gx{\Xi}         \def\Gy{\Psi}

%%Macros for calligraphic letters.
\def\CS{{\mathcal S}}   \def\CM{{\mathcal M}}   \def\CN{{\mathcal N}}
\def\CR{{\mathcal R}}   \def\CO{{\mathcal O}}   \def\CP{{\mathcal P}}
\def\CA{{\mathcal A}}   \def\CB{{\mathcal B}}   \def\CC{{\mathcal C}}
\def\CD{{\mathcal D}}   \def\CE{{\mathcal E}}   \def\CF{{\mathcal F}}
\def\CG{{\mathcal G}}   \def\CH{{\mathcal H}}   \def\CI{{\mathcal I}}
\def\CJ{{\mathcal J}}   \def\CF{{\mathcal F}}   \def\CL{{\mathcal L}}
\def\CT{{\mathcal T}}   \def\CU{{\mathcal U}}   \def\CV{{\mathcal V}}
\def\CZ{{\mathcal Z}}   \def\CX{{\mathcal X}}   \def\CY{{\mathcal Y}}
\def\CW{{\mathcal W}}
%%%%%
%%Macros for 'blackboard' letters (See (27) for display.)
\def\BBA {\mathbb A}   \def\BBb {\mathbb B}    \def\BBC {\mathbb C}
\def\BBD {\mathbb D}   \def\BBE {\mathbb E}    \def\BBF {\mathbb F}
\def\BBG {\mathbb G}   \def\BBH {\mathbb H}    \def\BBI {\mathbb I}
\def\BBJ {\mathbb J}   \def\BBF {\mathbb F}    \def\BBL {\mathbb L}
\def\BBM {\mathbb M}   \def\BBN {\mathbb N}    \def\BBO {\mathbb O}
\def\BBP {\mathbb P}   \def\BBR {\mathbb R}    \def\BBS {\mathbb S}
\def\BBT {\mathbb T}   \def\BBU {\mathbb U}    \def\BBV {\mathbb V}
\def\BBW {\mathbb W}   \def\BBX {\mathbb X}    \def\BBY {\mathbb Y}
\def\BBZ {\mathbb Z}

%%Macros for Ghotic (Fraktur) letters.
\def\GTA {\mathfrak A}   \def\GTB {\mathfrak B}    \def\GTC {\mathfrak C}
\def\GTD {\mathfrak D}   \def\GTE {\mathfrak E}    \def\GTF {\mathfrak F}
\def\GTG {\mathfrak G}   \def\GTH {\mathfrak H}    \def\GTI {\mathfrak I}
\def\GTJ {\mathfrak J}   \def\GTF {\mathfrak F}    \def\GTL {\mathfrak L}
\def\GTM {\mathfrak M}   \def\GTN {\mathfrak N}    \def\GTO {\mathfrak O}
\def\GTP {\mathfrak P}   \def\GTR {\mathfrak R}    \def\GTS {\mathfrak S}
\def\GTT {\mathfrak T}   \def\GTU {\mathfrak U}    \def\GTV {\mathfrak V}
\def\GTW {\mathfrak W}   \def\GTX {\mathfrak X}    \def\GTY {\mathfrak Y}
\def\GTZ {\mathfrak Z}   \def\GTQ {\mathfrak Q}

\font\Sym= msam10 % special symbols
\def\SYM#1{\hbox{\Sym #1}}
\newcommand{\bdw}{\prt\Gw\xspace}
%%%%%%%%%%%%%%%%%%%%%%%%%%%%%%%
\begin{abstract}
We give a global bilateral estimate on the maximal solution $\bar
u_F$ of $\,\prt_tu-\Gd u+u^q=0$ in $\BBR^N\ti (0,\infty)$, $q>1$,
$N\geq 1$, which vanishes at $t=0 $ on the complement of a closed
subset $F\subset \BBR^N$. This estimate is expressed by a Wiener
test involving the Bessel capacity $C_{2/q,q'}$. We  deduce from
this estimate that $\bar {u}_F$ is $\gs$-moderate in Dynkin's
sense.
\end{abstract}
\maketitle
\vskip -10mm
\noindent\underline{\phantom{----------------------------------------------------------------------------------------------------}}
\begin{center}{\sc\bf Repr\'esentation capacitaire des solutions positives\\
d'\'equations
  paraboliques semi-lin\'eaires}
\end{center}
\begin{quotation}
{\fnz {\sc{R\'esum\'e}.}
%\begin {center} {\bf R\'esum\'e}
%\end {center}
%{\small
%
\setcounter{equation}{0} Nous donnons une estimation bilat\'erale
pr\'ecise de la solution maximale $\bar {u}_F$ de $\prt_tu-\Gd
u+u^q=0$ dans $\BBR^N\ti (0,\infty)$, $q>1$, $N\geq 1$, qui
s'annulle en $t=0 $ sur le compl\'ementaire d'un sous-ensemble
ferm\'e  $F\subset \BBR^N$. Cette estimation s'exprime par un test
de Wiener impliquant la capacit\'e de Bessel $C_{2/q,q'}$. Nous
d\'eduisons de cette estimation que $\bar {u}_F$ est
$\gs$-moder\'ee au sens de Dynkin.}
\end {quotation}
\underline{\phantom{----------------------------------------------------------------------------------------------------}}

\vskip 2mm
%%%%%%%%%%%%%%%%%%%%%%%%%%%%%%%%%%%%%%%%%%%%%
%%%%%%%%%%%%%%%%%%%%%%%%%%%%%%%%%%%%%%%%%%%%%%%
%%%%FRENCHVERSION%%%%%%%%%%%%%%%%%%%%%%%%%%%%%%%%%%%%%%%%%%%%%%%%%%%%%%%%%%%%%%%%%%%%%%%%%%%%%
%\nind{\bf{Version fran\c caise abr\'eg\'ee}}\medskip
\section*{Version fran\c caise abr\'eg\'ee}
Soit $q>1$. Si $u$ est une solution positive de
\begin {equation}\label {mainf}
\prt_tu-\Gd u+\abs u^{q-1}u=0
\end {equation}
dans $\BBR^N\ti (0,\infty)$, nous avons d\'emontr\'e dans \cite
{MV1} qu'elle admet une trace initiale, not\'ee $Tr(u)$, dans la
classe des mesures de Borel positives et r\'eguli\`eres, mais pas
n\'ecessairement localement born\'ees. Si $F$ est un sous-ensemble
ferm\'e de $\BBR^N$, nous d\'esignons par $\bar {u}_F$ la solution
maximale de (\ref{mainf}) dont le support de la trace initiale est
inclus dans $F$. Si $1<q<q_c:=(N+2)/N$, il est montr\'e dans \cite
{BPT} que les in\'egalit\'es suivantes sont v\'erifi\'ees
$$
t^{-1/(q-1)}f(\abs {x-a}/\sqrt t)\leq  \bar {u}_F(x,t)\leq
\left((q-1)t\right)^{-1/(q-1)}\forevery a\in F,
$$
o\`u $f$ est l'unique fonction positive v\'erifiant
$$
\Gd f+\myfrac {1}{2}y\cdot Df+\myfrac {1}{q-1}f-\abs f^{q-1}f=0\; \mbox { dans
}\BBR^N\; \mbox { et } \lim_{\abs y\to\infty}\abs y^{2/(q-1)}f(y)=0.$$
Ces in\'egalit\'es jouent un r\^ole fondamental dans la d\'emonstration du
caract\`ere biunivoque de la correspondance, par l'op\'erateur de trace initiale, 
entre l'ensemble des solutions positives de (\ref {mainf}) et l'ensemble des mesures
de Borel
positives r\'eguli\`eres. Quand $q\geq q_c$ la fonction $f$ est identiquement nulle
car les singularit\'es isol\'ees
de (\ref {mainf}) sont \'eliminables  \cite {BF}.\medskip

\nind{\bf D\'efinition.} {\it Soit $N\geq 1$, $q\geq q_c$ et $F$ un sous-ensemble ferm\'e
de $\BBR^N$. On d\'efinit le potentiel $(2/q,q')$-capacitaire $W_{F}$ de $F$ par
\begin{equation}\label {F1}
W_{F}(x,t)=t^{-1/(q-1)}
\sum_{n=0}^\infty(n+1)^{N/2-1/(q-1)}e^{-n/4}
C_{2/q,q'}\left(\myfrac {F_{n}}{\sqrt {(n+1)t}}\right),
\end {equation}
 $\forall (x,t)\in\BBR^N\ti[0,\infty)$, o\`u $F_n=F_n(x,t)=\left\{y\in
F:\sqrt{nt}\leq\abs {x-y}< \sqrt {(n+1)t}\right\}$.}\medskip

 Notre r\'esultat principal est l'estimation bilat\'erale.\medskip

\nind{\bf Th\'eor\`eme 1.} {\it  Soit $N\geq 1$, $q\geq q_c$. Il existe deux constantes
$C_{1}>C_{2}>0$, ne d\'ependant que de $N$ et $q$, telles que pour tout
sous-ensemble ferm\'e $F$ de $\BBR^N$
\begin{equation}\label {F1'}
C_{2}W_{F}(x,t)\leq \bar {u}_{F}(x,t)\leq C_{1}W_{F}(x,t)\forevery
(x,t)\in \BBR^N\ti (0,\infty).
\end {equation}
}
\indent Si $\gm\in \mathfrak M^b(\BBR^N) $ (l'espace des mesures de Radon born\'ees
dans $\BBR^N$) appartient \`a $W^{-2/q,q}(\BBR^N)$, il existe une unique solution
$u={u}_\gm$ de (\ref {mainf}) dont la trace initiale est $\gm$ \cite {BP}.
On d\'efinit alors
\begin{equation}\label {F3}
\underline {u}_F=
\sup\{{u}_\gm:\gm\in \mathfrak M_+^b(\BBR^N)\cap W^{-2/q,q}(\BBR^N):\gm(F^c)=0  \}.
\end {equation}
Cette solution est $\gs$-mod\'er\'ee au sens de Dynkin, c'est \`a dire qu'il existe
une suite croissante de mesures positives $\gm_n\in \mathfrak M_+^b(\BBR^N)\cap
W^{-2/q,q}(\BBR^N)$ telles que ${u}_{\gm_n}\uparrow \underline {u}_F$.\medskip

La clef de la d\'emonstration du Th\'eor\`eme 1 est l'estimation inf\'erieure de
$\underline u_{F}$.\medskip

\nind{\bf Th\'eor\`eme 2.}  {\it Soit $N\geq 1$, $q\geq q_c$.  Il existe
une constante $C=C(N,q)>0$ telle que pour tout sous-ensemble
ferm\'e $F$ de $\BBR^N$,
\begin{equation}\label {F4} \underline u_{F}(x,t)\geq C W_{F}(x,t)\forevery
(x,t)\in \BBR^N\ti (0,\infty).
\end {equation}
}

Une cons\'equence importante des estimations pr\'ec\'edentes est la suivante:\medskip

\nind{\bf Th\'eor\`eme 3.}  {\it Soit $N\geq 1$, $q>1$.
Pour tout sous-ensemble ferm\'e $F$ de $\BBR^N$, $\bar {u}_{F}=\underline {u}_{F}$.}

%%%%%%%%%%%%%%%%%%%%%%%%%%%%%%%%%%%%%%%%%%%%%
%%%%%%%%%%%%%%%%%%%%%%%%%%%%%%%%%%%%%%%%%%%%%%%
%%%%ENGLISH VERSION%%%%%%%%%%%%%%%%%%%%%%%%%%%%%%%%
%%%%%%%%%%%%%%%%%%%%%%%%%%%%%%%%%%%%%%%%%%%%%
\setcounter{equation}{0}
%\nind{\large\bf{Main results}}\\ [-2mm]
\section{Main results}
Let  $u$ be a nonnegative solution of
\begin {equation}\label {main}
\prt_tu-\Gd u+\abs u^{q-1}u=0, \qq q>1
\end {equation}
in $\BBR^N\ti (0,\infty)$. It was proved in \cite {MV1} that $u$
admits an initial trace, denoted by $Tr(u)$, in the class of outer
regular positive Borel measures, not necessarily locally bounded.
If $F$ is a closed subset of $\BBR^N$, we denote by $\bar {u}_F$
the maximal solution of (\ref {main})  which belongs to $C(F^c\ti
[0,\infty))$ and vanishes on $F^c\ti \{0\}$, where
$F^c:=\BBR^N\sms F$.

 If $1<q<(N+2)/N$, the following inequalities are verified
\begin {equation}\label {main1}
t^{-1/(q-1)}f(\abs {x-a}/\sqrt t)\leq  \bar {u}_F(x,t)\leq
\left((q-1)t\right)^{-1/(q-1)}\forevery a\in F,
\end {equation}
where $f$ is the unique positive solution \cite {BPT} of
\begin {equation}\label {vss}
\Gd f+\myfrac {1}{2}y\cdot Df+\myfrac {1}{q-1}f-\abs f^{q-1}f=0\;
\mbox { in }\BBR^N\; \mbox { s. t. } \lim_{\abs y\to\infty}\abs
y^{2/(q-1)}f(y)=0.
\end {equation}
These inequalities play a fundamental role in proving that (in the subcritical case)
any positive solution is uniquely
determined by its initial trace \cite{MV1}.\medskip

\noindent{\bf Definition 1.} {\it Let
\begin{equation}\label {slic}\BA {l}
B_{\sqrt{nt}}(x)=\{y\in\BBR^N:\abs {x-y}<\sqrt{nt}\}\\[2mm]
T_n(x,t)=\bar B_{\sqrt{(n+1)t}}(x)\setminus B_{\sqrt{nt}}(x).
\EA\end {equation} For every $q\geq q_c$ and every closed subset
$F\subset \BBR^N$, we define the $(2/q,q')$-capacitary potential
$W_{F}$ of $F$ by
\begin{equation}\label {E1}
W_{F}(x,t)=t^{-1/(q-1)}
\sum_{n=0}^\infty(n+1)^{N/2-1/(q-1)}e^{-n/4}
C_{2/q,q'}\Big(\myfrac {F\cap T_n(x,t)}{\sqrt {(n+1)t}}\Big),
\end {equation}
 $\forall (x,t)\in\BBR^N\ti[0,\infty)$.
}\medskip

Our main result is the following bilateral estimate
\bth {Bila} Let $N\geq 1$, $q\geq q_c$. There exist two positive constants $C_{1}$ and
$C_{2}$, depending on $N$ and $q$ such that for any closed subset $F$ of $\BBR^N$,
\begin{equation}\label {E1'}
C_{2}W_{F}(x,t)\leq \bar {u}_{F}(x,t)\leq C_{1}W_{F}(x,t)\forevery
(x,t)\in \BBR^N\ti (0,\infty).
\end {equation}
\es

\nind\Remark It is important to notice that although $W_F$ is not a solution of
(\ref{main}),
it is invariant with respect to the similarity transformation associated with this
equation:
\begin{equation}\label {T1}
k^{1/(q-1)}W_{F}(\sqrt kx,kt)=W_{F/\sqrt k}(x,t)\forevery
(x,t)\in \BBR^N\ti (0,\infty),\;\forall k>0.
\end {equation}
Clearly $\bar {u}_{F}$ is also self similar, i.e., invariant with
respect to the above transformation.

%If $\gm\in \mathfrak M^b(\BBR^N) $ (the space of bounded Radon measures in
%$\BBR^N$) belonging to $W^{-2/q,q}(\BBR^N)$,
If $\gm$ is a bounded Borel measure such that $|\gm|\in W^{-2/q,q}(\BBR^N)$ then,
there exists a unique solution
$u={u}_\gm$ of (\ref {main}) with initial trace $\gm$.
We define
$$
\underline {u}_F=
\sup\{{u}_\gm:\gm\in  W_+^{-2/q,q}(\BBR^N),\;\gm(F^c)=0 \}.
$$
This solution is  $\gs$-moderate in the sense of Dynkin, which means that there
exists an increasing sequence of
positive measures $\gm_n\in  W_+^{-2/q,q}(\BBR^N)$ such that ${u}_{\gm_n}\uparrow
\underline {u}_F$. \medskip
 Clearly $\underline{u}_F\leq \bar u_F$. Therefore the next result implies the lower
estimate in
\rth {Bila}.
\bth {Dow} Let $N\geq 1$, $q\geq q_c$. There exists a positive $C=C(N,q)$,  such that for
any closed subset $F$ of $\BBR^N$,
\begin{equation}\label {E2}
\underline u_{F}(x,t)\geq C W_{F}(x,t)\forevery
(x,t)\in \BBR^N\ti (0,\infty).
\end {equation}
\es

As a consequence of Theorems \ref{t:Bila} and \ref{t:Dow} we find that $
\bar{u}_F\leq c \underline{u}_F$.
 Using this fact we obtain the following result (already known \cite {MV1} in the case $1<q<q_c$),
\bth {moder}
Let  $N\geq 1$, $q>1$. For any closed subset   $F$ of
$\BBR^N$ one has $ \bar {u}_F= \underline {u}_F$. In particular
$\bar {u}_F$ is $\gs$-moderate. \es
\medskip

%%%%%%%%%%%%%%%%%%%%%%%%%%%%%%%%%%%%%%%%%%%%%
%%%%%%%%%%%%%%%%%%%%%%%%%%%%%%%%%%%%%%%%%%%%%%%
%%%%ABRIDGED PROOF%%%%%%%%%%%%%%%%%%%%%%%%%%%%%%%%%%
%%%%%%%%%%%%%%%%%%%%%%%%%%%%%%%%%%%%%%%%%%%

\section{ Proof of the upper estimate in \rth{Bila}} In the sequel we denote by $c$
a positive constant which  depends only on $N$ and $q$;
its value may change from one occurrence to another. Without loss of generality, we
assume that $F$ is compact.
Denote $B_r(x)=\{y:\abs {x-y}<r\}$ and $B_r=B_r(0)$.
\par Let $r>0$ be a positive number such that $F\subset B_r$. We start by deriving
an upper estimate depending on $r$.
Let $\gr>0$ be a positive number,  to be later determined as a function of $r,t$.
Let $\eta\in C^\infty_0(B_{r+\gr})$ be such that $\eta=1$ on $F$ and $0\leq\eta\leq 1$.
Put $\eta^*=1-\eta$ and choose
$\gz:=\left(e^{t\Gd}[\eta^*]\right)^{2q'}$ as a test function. Then
\begin{equation}\label{int-main}
\int_0^1 \int_{\BBR^N} u(\prt_t -\Gd)\gz dxdt + \int_0^1\int_{\BBR^N} u^q\gz
dxdt=-\int_{\BBR^N}u(x,1)dx.
\end{equation}
A straightforward computation yields
$$\int_{Q_{r+\gr}} u^q dxdt + \int_{\BBR^N}u(x,1) dx\leq
\int_0^1\int_{\BBR^N}\left(R(\eta)\right)^{q'}\,dx\,dt,$$
$$R(\eta):=\abs {De^{t\Gd}[\eta]}^2+\abs {\prt_te^{t\Gd}[\eta]+\Gd e^{t\Gd}[\eta]},\q
Q_{r}:=\{(x,t):t>0,\;\abs x^2+t\geq r^2\}.$$
Using interpolation inequalities \cite {Tr} one obtains,
\begin{equation}\label{up1}
  \int_0^1\int_{\BBR^N}\left(R(\eta)\right)^{q'}\,dx\,dt
\leq C_1\norm\eta_{W^{2/q,q'}}^{q'}.
\end{equation}
This implies
$$\int_{\BBR^N}u(x,1)\,dx+\int_{Q_{r+\gr}}u^q\,dx\,dt
\leq cC_{2/q,q'}^{B_{r+\gr}}(F).
$$
Further it can be shown that, for $0<s<1$,
\begin{equation}\label{up2}
\int_{\BBR^N}u(x,1)\,dx+\int_{s}^1\int_{\BBR^N}u^q\,dx\,dt
=\int_{\BBR^N}u(x,s)\,dx.
\end{equation}
Clearly $u(x,t+s)<w_s(x,t):=e^{t\Gd}[u(\cdot,s)]$ for $t>0$. Therefore, by
(\ref{up1}) and (\ref{up2}),
\begin{equation}\label{up3}
 u(x,(r+2\gr)^2)\leq \myfrac{c}{(\gr^2+r\gr)^{N/2}}C_{2/q,q'}^{B_{r+\gr}}(F).
\end{equation}
%Put $Q_{r,\gr}^*=\{(x,t):\abs x> r+\gr,\; 0\leq t< (r+2\gr)^2\}$ and let
Let $v$ be the solution of the initial-boundary value problem
$$\BAL \prt_tv-\Gd v=0 &\txt{in} Q^*_{r,\gr}:=(\bar B_{r+\gr})^c\ti (0,(r+2\gr)^2), \\
v(x,0)=0 \forevery x\in B^c_{r+\gr},\q &v=u \forevery (x,t)\in \prt B_{r+\gr}\ti
(0,(r+2\gr)^2).
\EAL$$
Using the fact that $u<v$ in $Q^*_{r,\gr}$ and (\ref{up3}) we obtain,

\blemma{uplem} If $F\subset B_r$, there exists $c>0$ such that
\begin {equation}\label {upest1}
\bar {u}_F(x,t)\leq c\left(1+\myfrac{r}{\gr}\right)^{N/2}\myfrac
{e^{-(\abs x-r-3\gr)^2/4t}}{t^{N/2}}C_{2/q,q'}^{B_{r+\gr}}(F),
\end {equation}
for any $(x,t)\in \BBR^N\ti[(r+3\gr)^2,\infty)$.
\es
The upper estimate in (\ref{E1'}) is obtained by slicing $F$, relative to a given
point $(x,t)\in\BBR^N\ti (0,\infty)$,
in such a way that each slice
satisfies the assumption of \rlemma{uplem} for an appropriate value of $r$ depending
on the point.
Put
$$T_n(x,t)=B_{\sqrt {(n+1)t}}(x)\setminus B_{\sqrt {nt}}(x),\;F_n(x,t)=F\cap
T_n(x,t)\forevery n\in\BBN.
$$
Since a sum of positive solutions of (\ref {main}) is a super-solution,
$$\bar u_F\leq\sum_{n=0}^\infty \bar u_{F_n(x,t)}.
$$
Using this fact and \rlemma{uplem} we show that
%t^{-1/(q-1)}\sum_{n=0}^\infty(n+1)^{N/2-1/(q-1)}
\begin{equation}\label{up-n}
  \bar u_F(x,t)\leq cW_F(x,t)
  %\myfrac {c}{t^{1/(q-1)}}
%\sum_{n=0}^\infty(n+1)^{N/2-1/(q-1)}e^{-n/4}C_{2/q,q'}\left(\myfrac
%{F_{n}(x,t)}{\sqrt {(n+1)t}}\right)
\end{equation}
for every $(x,t)\in \BBR^N\ti (0,\infty)$.
In view of the fact that both sides of this inequality are  invariant with respect
to the similarity transformation
(\ref{T1}), it is sufficient to prove it in the case $(x,t)=(0,1)$. We denote
$F_n=F_n(0,1)$ and
$T_n=T_n(0,1)$.
\par If $N=1$, each of the sets $F_n$ satisfies the condition of \rlemma{uplem} with
$r=\sqrt n$.
But this is not the case when $N\geq 2$. Therefore, if $N\geq 2$, a secondary
slicing is needed.
\par  For every $n\in \BBN$ there exists  a set of points
$\Gth_{n}=\{a_{j,n}\}_{j=1}^{J_n}$
 on the sphere $\abs{y}=(\sqrt{n+1}+\sqrt n)/2$ such that
$$
\abs {a_{n,j}-a_{n,k}}\geq 1/\sqrt{2(n+1)} \txt{for} j\neq k, \q
T_{n}\subset\bigcup_{1\leq j\leq J_n}B_{\sqrt {1/(n+1)}}(a_{n,j}).$$
Clearly $J_n\leq (\sqrt{2}(n+1))^{N-1}< (4n)^{N-1}$.
If $F_{n,j}:=F_{n}\cap B_{\sqrt {1/(n+1)}}(a_{n,j}) $,
\begin{equation}\label{well}
\bar {u}_{F}(0,1)\leq \sum_{n=0}^\infty\sum_{1\leq j\leq J_n}\bar {u}_{F_{n,j}}(0,1).
\end {equation}
It is not difficult to verify that
$$
C_{2/q,q'}^{B_{2/\sqrt {n+1}}(a_{n,j})}(F_{n,j})
\approx (n+1)^{N/2-1/(q-1)}
C_{2/q,q'}\left(\sqrt {n+1}\,F_{n,j}\right)
$$
where the capacity on the right hand side (resp. left hand side) is the Bessel
capacity relative
to $\BBR^N$ (resp. relative to $B_{2/\sqrt {n+1}}(a_{n,j})$).
The symbol $\approx$ \ stands for two-sided inequalities with constants $c=c(N,q)$.
Further,
 the quasi-additivity of Bessel capacities \cite {AB} implies
\begin{equation}\label{well2}
\sum_{1\leq j\leq J_n}C_{2/q,q'}(\sqrt {n+1}\,F_{n,j})
\leq c(N,q)\,C_{2/q,q'}(\sqrt {n+1}\,F_{n}).
\end {equation}
Each set $F_{n,j}$ satisfies the condition of \rlemma{uplem} with $r=\sqrt n$.
Therefore,
estimating $\bar {u}_{F_{n,j}}(0,1)$ as in (\ref{upest1}) and using
(\ref{well}) and (\ref{well2}) we obtain (\ref {up-n}) for $(x,t)=(0,1)$.
%%%%%%DOWN***PROOF%%%%%%%%%%%%%%%%%%%%%%%%%%%%

\section{Proof of \rth{Dow}}
 As in the elliptic case \cite{MV2},  for each $(x,t)\in \BBR^N\ti (0,\infty)$, we
construct
 a measure $\gm=\gm_{x,t}\in W_+^{-2/q,q}(\BBR^N)$, concentrated on $F$, such that
\begin{equation}\label {sum2}
 {u}_\gm(x,t)\geq cW_F(x,t).
% \myfrac {c}{t^{1/(q-1)}}
%\sum_{n=0}^\infty(n+1)^{N/2-1/(q-1)}e^{-n/4}C_{2/q,q'}\left(\myfrac {F_n}{\sqrt
%{(n+1)t}}\right).
 \end{equation}
Since $\underline u_F> u_{\gm}$, (\ref{sum2}) implies (\ref{E2}).
%As before, we may assume that $(x,t)=(0,1)$.
\par For every measure $\gm\in W_+^{-2/q,q}(\BBR^N)$,
 $0\leq {u}_\gm\leq e^{t\Gd}[\gm]$. Therefore
\begin{equation}\label{Ineq1}
  {u}_\gm\geq e^{t\Gd}[\gm]-\int_0^te^{(t-s)\Gd}[(e^{s\Gd}[\gm])^q]ds.
\end{equation}
Let $\gn_n$ be the capacitary measure of $F_{n}(x,t)/\sqrt {t(n+1)}$ (see \cite{AH})
and define the measure
$\mu_n$ by
$$\gm_n(A)=(t(n+1))^{N/2-1/(q-1)}\gn_n({A}/{\sqrt {t(n+1)}}),$$
for every Borel set $A$. Thus $\mu_n$ is concentrated on $F_n$. Finally put
$\mu:=\mu_{x,t}=\sum_n\gm_n$. With this choice of $\mu$ it is not difficult to show
that
%obtain the following estimate from below for the linear term in (\ref{Ineq1}):
%$$\gn_n( {F_n}/{\sqrt {t(n+1)}})
%=C_{2/q,q'}({F_n}/{\sqrt {t(n+1)}})=\norm{\gn_n}^q_{W^{-2/q,q}}$$
\begin{equation}\label{est4}
e^{t\Gd}[\gm](x,t)\geq
\dfrac {1}{(4\gp t)^{N/2}}\sum_{n=0}^{\infty}
(\sqrt{(n+1)t})^{N-2/(q-1)}e^{-(n+1)/4} C_{2/q,q'}\Big(\dfrac
{F_{n}(x,t)}{\sqrt {(n+1)t}}\Big).
\end{equation}
We have to derive a corresponding upper estimate for
the nonlinear term in (\ref{Ineq1}).
 To this end we employ a partitioning $\{\CT_n:n\in \BBZ\}$ of $\BBR^N\ti (0,t)$
defined by,
 $$\CT_n=\begin{cases}
 \{(y,s): {tn}\leq  {\abs {x-y}^2+t-s}\leq {t(n+1)},\;0<s<t\}, &\text{ if }
n\in\BBN_*,\\
  \{(y,s): t \ga^{-n}\leq \abs {x-y}^2+t-s\leq t \ga^{-n-1},\;0<s<t\},
  &\text{ if } n\leq 0,
 \end{cases}$$
 where $\ga\in (0,1)$ must be appropriately chosen. If
$G(\gx,\gt):=(4\gp|\gt|)^{-N/2}\exp (-\abs{\gx}^{2}/4\gt))$ then
  $$
  \int_{0}^{t}e^{-(t-s)\Gd}(e^{s\Gd}[\gm])^q ds=
  C\sum_{p\in\BBZ}\iint_{\CT_p}G(y-x,s-t)
  \Big(\sum_{n=0}^\infty \int_{\BBR^N} G(z-y,s)d\gm_n(z)\Big)^q dyds.
  $$
We denote
 $$J_1=\sum_{p\in\BBZ}\iint_{\CT_p} G(y-x,s-t)
  \Big(\sum_{n=0}^{p+2}
 \int_{\BBR^N} G(z-y,s){s^{N/2}}d\gm_n(z)\Big)^q dyds,
 $$
  $$J_2=\sum_{p\in\BBZ}\iint_{\CT_p} G(x-y,s-t)
  \Big(\!\!\sum_{n=p+3}^{\infty}
  \int_{\BBR^N} G(z-y,s){s^{N/2}}d\gm_n(z)\Big)^q dyds.
 $$
\begin{comment}
 %%%%%%%%%%%%%%%%%%%%%
  $$
  \int_{0}^{t}e^{-(t-s)\Gd}(e^{s\Gd}[\gm])^q ds=
  C\sum_{p\in\BBZ}\iint_{\CT_p}\dfrac{e^{\abs {x-y}^2/4(t-s)}}{(t-s)^{N/2}}
  \Big(\sum_{n=0}^\infty \int_{\BBR^N}\dfrac{e^{-\abs
{y-z}^2/4s}}{s^{N/2}}d\gm_n(z)\Big)^q ds.
 $$
We denote
 $$J_1=\sum_{p\in\BBZ}\iint_{\CT_p}
\dfrac {e^{\abs {x-y}^2/4(t-s)}}{(t-s)^{N/2}}
  \Big(\sum_{n=0}^{p+2}
 \dfrac {e^{-\abs {y-z}^2/4s}}{s^{N/2}}d\gm_n(z)\Big)^qds,
 $$
  $$J_2=\sum_{p\in\BBZ}\iint_{\CT_p}
\dfrac {e^{\abs {x-y}^2/4(t-s)}}{(t-s)^{N/2}}
  \Big(\!\!\sum_{n=p+3}^{\infty}
 \dfrac {e^{-\abs {y-z}^2/4s}}{s^{N/2}}d\gm_n(z)\Big)^qds.
 $$
%%%%%%%%%%%%%%%%%%%%%%%%%%%%%%%
\end{comment}
 $J_1$ is estimated  using H\"older's inequality and the following inequalities
\cite {Tr}
 \begin {equation}\label{intest}
 \myfrac {1}{C}\norm{\gl}_{W^{-2/q,q}}
 \leq \norm {e^{t\Gd}[\gl]}_{L^q(\BBR^N\ti (0,1))}
 \leq {C\norm{\gl}_{W^{-2/q,q}}},
\end {equation}
 valid for any $\gl\in W^{-2/q,q}$. The estimate of $J_2$, in the case $N>1$,
requires a more
 delicate argument using the secondary slicing introduced in the previous section
and the
 quasi-additivity of Bessel capacities. For a suitable choice of $\ga$ one obtains
 \begin {equation}\label {est5}
 J_1+J_2\leq
\dfrac {C}{t^{N/2}}\sum_{n=0}^{\infty}
(\sqrt{(n+1)t})^{N-2/(q-1)}e^{-(n+1)/4} C_{2/q,q'}\Big(\dfrac{F_{n}}{\sqrt
{(n+1)t}}\Big).
\end {equation}
The inequalities (\ref{Ineq1}), (\ref{est4}) and (\ref{est5}) imply (\ref{sum2}), with
 $\mu$ replaced by $\ge\mu$,  provided that $\ge>0$ is
sufficiently small, depending on $q$ and the constants in (\ref{est4}) and
(\ref{est5}).
This in turn implies (\ref{E2}).
\par It follows that $$\underline {u}_F(x,t)\leq \bar {u}_F(x,t)\leq c\underline
{u}_F(x,t).$$
By the uniqueness argument used in \cite {MV1} we conclude that $\bar {u}_F(x,t)=
\underline {u}_F(x,t).$

\vskip 3mm
\noindent{\bf Acknowledgment.}\hskip 2mm The authors are partially sponsored by an
E.C. grant through
the RTN program 'Front Singularities' HPRN-CT-2002-00274.

%%%END DOCUMENT%%%%%%%%%%%%%%%

\begin{thebibliography}{99}

\bibitem{AH} Adams D. R. and Hedberg L. I., Function Spaces and Potential Theory,
Grundlehren  Math. Wissen. {\bf 314}, Springer (1996).

\bibitem {AB} Aikawa H. and Borichev A. A., {\em Quasiadditivity and measure
property of capacity and the tangential boundary behaviour of harmonic functions},
Trans. Amer. Math. Soc. {\bf 348}, 1013-1030 (1996).
\bibitem{BP} Baras P. and Pierre M., Probl\'{e}mes paraboliques semi-lin\'{e}aires
avec donn\'{e}es mesures, {\it Applicable Anal.} {\bf 18} (1984), 111-149.
\bibitem {BF} Brezis H. and Friedman A., {\em Nonlinear parabolic equations involving
measures as initial data}, J. Math. Pures Appl. {\bf 62}, 73-97 (1983).

\bibitem {BPT} Brezis H., Peletier  L. A. and Terman D., {\em A very singular equation
of the heat equation with absorption}, Arch. Rat. Mech. Anal. {\bf 95}, 185-209
(1986).

\bibitem{MV1} Marcus M. and V\'{e}ron L., {\em  The initial trace of positive
solutions of semilinear parabolic equations}, Comm. in P.D.E. {\bf 24}, 1445-1499
(1999).

\bibitem{MV2} Marcus M. and V\'{e}ron L., {\em
Capacitary estimates of positive solutions of
semilinear elliptic equations with absorption}, J. Eur. Math. Soc. {\bf 6}, 483-527
(2004).

\bibitem{Tr} Triebel H., Interpolation Theory, Function Spaces, Differential
Operators, North-Holland Amsterdam (1978).
\end{thebibliography}
\end {document}